\documentclass[10pt]{article}

\usepackage[latin1]{inputenc}
\usepackage{amsmath}
\usepackage{amsfonts,amssymb,latexsym}
\usepackage[francais]{babel}
\usepackage{epsfig,epsf}
\newcounter{fig}
\def\figcaption #1
  {\addtocounter{fig}{1}
   \begin{center}
   {\sc Fig.} \arabic{fig}. {\it #1} 
    \end{center}
    \addtocontents{lof}{\protect \contentsline {figure}{\protect \numberline {\thefig}{\protect \ignorespaces \protect \pit #1 }}{\thepage}}
   }
\input{cyracc.def}

\newtheorem{theo}{Th\'eor\`eme}

\newtheorem{prop}{Proposition}

\newcommand{\eps}{\varepsilon}

\newcommand{\ioe}{\leqslant}
\newcommand{\soe}{\geqslant}
\newcommand{\vers}{\rightarrow}

\newcommand{\intc}{{\frac{1}{2i\pi}\int \! \!}}

\newcommand{\dlz}{{\frac{\zeta '}{\zeta}}}
\newcommand{\demi}{{\frac{1}{2}}}

\newcommand{\Scal}{{\mathcal S}}

\newcommand{\Nat}{{\mathbb N}}

\newcommand{\Real}{{\mathbb R}}
\newcommand{\Com}{{\mathbb C}}

\newcommand{\card}{{\rm card}}

\newcommand{\fin}{\hfill$\Box$}
\newcommand{\dem}{\noindent {\bf D\'emonstration\ }}
\newcommand{\fine}{\tag*{\mbox{$\Box$}}}

\title{Notes de lecture de l'article\\ {\og Partial sums of the Möbius function\fg}\\ de Kannan Soundararajan}
\author{Michel Balazard et Anne de Roton}

\begin{document}
\maketitle

Dans ce document, nous exposons la démonstration du théorème suivant.
\begin{theo}[Soundararajan \cite{S2008}](HR)

Pour tout $\eps >0$, on a
$$
M(N)=\sum_{n \ioe N} \mu(n) \ll_{\eps} \sqrt{N}\exp \bigl ( (\log N)^{1/2} (\log \log N)^{5/2+\eps}\bigr ).
$$  
\end{theo}

Les lettres \textit{(HR)} indiquent que l'on suppose l'hypothèse de Riemann vérifiée.

La méthode de démonstration a été inventée par Maier et Montgomery dans l'article \cite{MM2008}. Ils y démontrent que
$$
M(N) \ll \sqrt{N}\exp \bigl ( (\log N)^{39/61}\bigr )
$$
sous l'hypothèse de Riemann. Leur approche a été ensuite perfectionnée par Soundararajan (cf. \cite{S2008}), qui a obtenu l'estimation
$$
M(N) \ll \sqrt{N}\exp \bigl ( (\log N)^{1/2} (\log \log N)^{14}\bigr ),
$$
toujours sous l'hypothèse de Riemann. Dans ce qui suit, nous suivons la démonstration de Soundararajan dans \cite{S2008}, en y incorporant les quelques précisions qui permettent de remplacer $14$ par n'importe quel nombre $>5/2$.

Un mot sur la locution {\og $T$ assez grand\fg}, que nous emploierons librement. Elle signifie $T \soe T_0$, où $T_0$ est une constante absolue et effectivement calculable (mais certainement très grande) telle que

$\bullet$ toutes les quantités dont nous parlons sont bien définies ;

$\bullet$ toutes les inégalités que nous écrivons sont vérifiées.

Par exemple, pour $T$ assez grand on a
$$
\frac{\log \log \log V'}{\log \log V'} \ioe \frac{\log \log \log V}{\log \log V}\quad (V'\soe V \soe (\log \log T)^2).
$$

Nous remercions Kannan Soundararajan pour de nombreux éclaircissements concernant son article \cite{S2008}.

\newpage

\tableofcontents

\newpage

\section{Ordonnées $V$-typiques de taille $T$}\label{t36}

Nous allons évaluer $M(N)$ grâce à la formule de Perron en utilisant un contour d'intégration sur lequel les grandes valeurs de $|\zeta (z)|^{-1}$ seront aussi rares que possible. Pour quantifier cette rareté, Soundararajan a introduit la notion suivante.

On se donne un paramètre $\delta$, tel que $0 <\delta \ioe 1$. Soit $T$ assez grand\footnote{En termes absolus, indépendamment de $\delta$.} et $V$ tel que $(\log \log T)^2 \ioe V \ioe \log T/\log\log T$. Un nombre réel $t$ est appelé une \textbf{ordonnée $V$-typique de taille $T$} si

$\bullet$ $T \ioe t \ioe 2T$ ; \label{t59}

\textit{(i)} pour tout $\sigma \soe 1/2$, on a
$$
\Bigl \lvert \sum_{n \ioe x} \frac{\Lambda (n)}{n^{\sigma +it}\log n}\frac{\log (x/n)}{\log x} \Bigr \rvert \ioe 2V, \quad \text{où $x=T^{1/V}$} ;
$$

\textit{(ii)} tout sous-intervalle de $[t-1,t+1]$ de longueur $2\delta\pi V/\log T$ contient au plus $(1+\delta) V$ ordonnées de zéros de $\zeta$ ;

\textit{(iii)} tout sous-intervalle de $[t-1,t+1]$ de longueur $2\pi V/\bigl ((\log V)\log T\bigr )$ contient au plus $V$ ordonnées de zéros de $\zeta$.

Si $t\in [T,2T]$ ne vérifie pas l'une des assertions \textit{(i), (ii), (iii)}, on dira que $t$ est une \textbf{ordonnée $V$-atypique de taille $T$}.

La pertinence de cette définition\footnote{Pour être précis, il faudrait parler d'ordonnée $(\delta,V)$-typique de taille $T$. Dans ce qui suit, la référence à $\delta$ sera implicite.} quant à la taille de $|\zeta (s)|$ est fournie par l'énoncé suivant.
\begin{prop}\label{t40}(HR)
 Soit $T$ assez grand et $V$ tel que $(\log \log T)^2 \ioe V \ioe \log T/\log\log T$. Soit $t$ une ordonnée $V$-typique de taille $T$. On a
$$
\log |\zeta (\sigma +it)| \soe -V\log \Bigl( \frac{V/\log T}{\sigma-1/2} \Bigr )-2(1+\delta)V\log\log V +O(V\delta^{-2})\quad\text{si} \quad \demi<\sigma \ioe \demi +\frac{V}{\log T},
$$
et
$$
\log |\zeta (\sigma +it)| \soe O(V\delta^{-1}) \quad \text{si} \quad  \demi +\frac{V}{\log T}\ioe \sigma \ioe 2.
$$  
\end{prop}

Cette proposition découle des propositions \ref{t60} et \ref{t61}, démontrées ci-dessous après une étude préliminaire de la dérivée logarithmique de la fonction $\zeta$.

\section{Le logarithme de la fonction $\zeta$}

\subsection{Estimations de $\dlz$}
Notons $\psi=\Gamma'/\Gamma$ la dérivée logarithmique de la fonction $\Gamma$. On a (cf. \cite{D2000}, chapter 12, (8), (11))
\begin{equation}
  \label{t53}
\dlz(s)=-\frac{1}{s}-\frac{1}{s-1}+\demi \log \pi -\demi \psi(s/2) +\sum_{\rho}  \frac{1}{s-\rho},  
\end{equation}
où la somme porte sur les zéros non triviaux $\rho=\beta +i\gamma$ de la fonction $\zeta$, et est calculée par la formule $\sum_{\rho}=\lim_{T\vers \infty}\sum_{|\gamma| \ioe T}$.

On en déduit d'abord la proposition suivante.
\begin{prop}
  Si $T$ est assez grand, on a
$$
\Re \dlz (\sigma +it) =F(s)-\demi \log T +O(1) \quad (1/2 \ioe \sigma \ioe 2, \, T\ioe t \ioe 2T, \, \zeta(\sigma+it)\not =0),
$$
où
\begin{equation}
  \label{t50}
  F(s)=\sum_{\rho} \Re \frac{1}{s-\rho}=\sum_{\rho} \frac{(\sigma-\beta)}{(\sigma-\beta)^2+(\tau-\gamma)^2}.
\end{equation}
\end{prop}
\dem

Cela résulte de l'identité \eqref{t53} et de l'estimation
\begin{equation*}
\psi (s)=\log s +O(s^{-1}) \quad (\Re s >0).\fine  
\end{equation*}

Nous donnons maintenant une autre identité concernant la dérivée logarithmique de la fonction zêta.

\begin{prop}\label{t54}
  Soit $x \soe 1$ et $z \in \Com$. On suppose que $z$ n'est pas un pôle de $\dlz$. Alors
$$
\sum_{n \ioe x} \frac{\Lambda (n)}{n^z} \log (x/n) = -\dlz(z)\log x -\Bigl ( \dlz \Bigr )'(z) -\sum_{\rho} \frac{x^{\rho -z}}{(\rho -z)^2} +\frac{x^{1-z}}{(1 -z)^2} -\sum_{n\soe 1}\frac{x^{-2n-z}}{(z+2n)^2}.
$$
\end{prop}
\dem

On peut supposer $x>1$ (pour $x=1$, l'identité se ramène à la dérivée de \eqref{t53}).
 
On constate que les deux membres sont des fonctions holomorphes de $z$ dans $\Com$ privé des pôles de $\dlz$. Il suffit donc de démontrer l'égalité pour $z$ réel $>1$. La méthode exposée au chapter 12 de \cite{D2000} s'applique, plus facilement car les intégrales seront absolument convergentes. On utilise la formule de Perron suivante :
$$
\sum_{n \ioe x} \frac{\Lambda (n)}{n^z} \log (x/n) = \frac{1}{2\pi i}\int_{c-i\infty}^{c+i\infty}-\dlz (z+w) x^w \frac{dw}{w^2} \quad (c>0).
$$

On déplace la droite d'intégration vers la gauche\footnote{Avec les précautions d'usage pour passer entre les zéros de $\zeta$, cf. \cite{D2000}, p.108.}, à l'abscisse $c'$ telle que $z+c'=-2N-1$, $N$ entier positif. D'après le théorème des résidus, on obtient
\begin{multline*}
\sum_{n \ioe x} \frac{\Lambda (n)}{n^z} \log (x/n) =  -\dlz(z)\log x -\Bigl ( \dlz \Bigr )'(z) -\sum_{\rho} \frac{x^{\rho -z}}{(\rho -z)^2} +\\
\frac{x^{1-z}}{(1 -z)^2} -\sum_{n \ioe N}\frac{x^{-2n-z}}{(z+2n)^2}+
\frac{1}{2\pi i}\int_{c'-i\infty}^{c'+i\infty}-\dlz (z+w) x^w \frac{dw}{w^2}.  
\end{multline*}

La dernière intégrale tend vers $0$ quand $N$ tend vers l'infini, en vertu de l'estimation 
$$
\dlz (s) \ll \log (2|s|),
$$
valable uniformément dans le demi-plan $\Re s \ioe -1$, privé des disques de rayon $1/2$, centrés en $-2,-4, \dots$.\fin\\

\begin{prop}
  Soit

$\bullet$ $T \soe 1$ ;

$\bullet$ $z \in \Com$ tel que $\Re z \soe 0$, $T \ioe |\Im z | \ioe 2T$ et $z$ n'est pas un pôle de $\dlz$ ;

$\bullet$ $1 \ioe x \ioe T$.

 Alors
 \begin{equation}
   \label{t55}
 \sum_{n \ioe x} \frac{\Lambda (n)}{n^z} \log (x/n) = -\dlz(z)\log x -\Bigl ( \dlz \Bigr )'(z) -\sum_{\rho} \frac{x^{\rho -z}}{(\rho -z)^2} +O(T^{-1}).  
 \end{equation}
\end{prop}
\dem

On applique la proposition \ref{t54} et on vérifie que
$$
\Bigl \lvert\frac{x^{1-z}}{(1-z)^2} \Bigr \rvert \ioe \frac{x}{T^2},
$$
et
\begin{align*}
\Bigl \lvert \sum_{n \soe 1}\frac{x^{-2n-z}}{(z+2n)^2}  \Bigr \rvert &\ioe  \sum_{n \soe 1}\frac{1}{T^2+4n^2}\\
& \ioe \frac{1}{T}.\fine
\end{align*}

\subsection{Estimations de $\log |\zeta|$ : démonstration de la proposition \ref{t40}}\label{t49}

\begin{prop}(HR)\label{t56}
  Soit $T$ assez grand et $T\ioe t\ioe 2T$. On a uniformément pour $1/2 <\sigma \ioe 2$ et $2 \ioe x \ioe T$
  $$
\log |\zeta (\sigma +it)| \soe \Re \sum_{n\ioe x} \frac{\Lambda (n)}{n^{\sigma +it}\log n} \frac{\log (x/n)}{\log x}- \Bigl ( 1 +\frac{x^{\demi-\sigma}}{(\sigma-1/2)\log x}\Bigr ) \frac{F(\sigma+it)}{\log x} +O(1).    
  $$
\end{prop}
\dem

On intègre \eqref{t55} entre $z=\sigma +it$ et $z=2+it$ :
\begin{multline*}
 \sum_{n \ioe x} \Lambda (n) \log (x/n)\Bigl ( \frac{n^{-2-it}}{-\log n}- \frac{n^{-\sigma -it}}{-\log n}\Bigr ) = \bigl ( -\log \zeta (2+it) + \log \zeta (\sigma +it) \bigr )\log x - \dlz(2+it) + \dlz (\sigma +it)\\
 -\sum_{\rho}\int_{\sigma}^2 \frac{x^{\rho -u-it}}{(\rho -u-it)^2}du +O(T^{-1}),  
\end{multline*}
donc
$$
(\log x) \log \zeta (\sigma +it) =  \sum_{n\ioe x} \frac{\Lambda (n)}{n^{\sigma +it}\log n} \log (x/n)- \dlz (\sigma +it) +\sum_{\rho}\int_{\sigma}^2 \frac{x^{\rho -u-it}}{(\rho -u-it)^2}du +O(\log x).
$$

On en déduit
\begin{align*}
  \log |\zeta (\sigma +it)| &=  \Re \sum_{n\ioe x} \frac{\Lambda (n)}{n^{\sigma +it}\log n} \frac{\log (x/n)}{\log x}-(\log x)^{-1}\Re\dlz (\sigma +it) +(\log x)^{-1}\Re \sum_{\rho}\int_{\sigma}^2 \frac{x^{\rho -u-it}}{(\rho -u-it)^2}du +O(1)\\
&=\Re \sum_{n\ioe x} \frac{\Lambda (n)}{n^{\sigma +it}\log n} \frac{\log (x/n)}{\log x}-(\log x)^{-1}F(\sigma+it)+\frac{\log T}{2\log x} +\\
& \quad +(\log x)^{-1} \Re \sum_{\rho}\int_{\sigma}^2 \frac{x^{\rho -u-it}}{(\rho -u-it)^2}du +O(1).
\end{align*}

Pour majorer le module de la somme sur les zéros, nous utilisons l'hypothèse de Riemann :
\begin{align*}
\Bigl \lvert \sum_{\rho}\int_{\sigma}^2 \frac{x^{\rho -u-it}}{(\rho -u-it)^2}du \Bigr \rvert & \ioe  \sum_{\rho}\frac{1}{|\rho -\sigma-it|^2}\int_{\sigma}^2 x^{\demi-u}du \\
& \ioe \frac{x^{\demi-\sigma}}{\log x}\sum_{\rho}\frac{1}{|\rho -\sigma-it|^2}\\\
&= \frac{x^{\demi-\sigma}}{(\sigma-1/2)\log x}F(\sigma+it),
\end{align*}
d'où le résultat annoncé.\fin

\medskip

Nous utiliserons dans la suite l'estimation simple suivante.

\begin{prop}\label{t57}
  Pour $a>0$, $c>0$ et $N \in \Nat$, on a
$$
\sum_{n=0}^N\frac{a}{a^2+(cn)^2} \ioe \frac{1}{a} +\frac{\pi}{2c}.
$$
\end{prop}
\dem

On a
\begin{align*}
\sum_{n =0}^N\frac{a}{a^2+(cn)^2} &=  \frac{1}{a} + \sum_{0<n\ioe N}\frac{a}{a^2+(cn)^2}\\
& \ioe  \frac{1}{a} + \int_0^{\infty}\frac{a}{a^2+(ct)^2}dt\\
& = \frac{1}{a} +\frac{a}{c^2}\int_0^{\infty}\frac{1}{(a/c)^2+t^2}dt\\
&= \frac{1}{a} +\frac{\pi}{2c}.
\end{align*}

\begin{prop}\label{t60}(HR)
Soit 

$\bullet$ $T$ assez grand ;

$\bullet$  $V$ tel que $(\log \log T)^2 \ioe V \ioe \log T/\log\log T$ ;

$\bullet$ $t$ une ordonnée $V$-typique de taille $T$.

Alors 
$$
\log |\zeta (\sigma +it)| \soe O(V/\delta) \quad ( 1/2+V/\log T  \ioe \sigma \ioe 2).
$$  
\end{prop}
\dem

Posons $x=T^{1/V}$. On a $2 \ioe x \ioe T$ et
\begin{align*}
\frac{x^{\demi-\sigma}}{(\sigma-1/2)\log x} &\ioe \frac{\exp(-V\log x/\log T)}{V\log x/\log T} \quad  ( 1/2+V/\log T  \ioe \sigma)\\
&= e^{-1} \ioe 1.  
\end{align*}

La proposition \ref{t56} donne
\begin{align*}
\log |\zeta (\sigma +it)| &\soe  \Re \sum_{n\ioe x} \frac{\Lambda (n)}{n^{\sigma +it}\log n} \frac{\log (x/n)}{\log x}- \Bigl ( 1 +\frac{x^{\demi-\sigma}}{(\sigma-1/2)\log x}\Bigr ) \frac{F(\sigma+it)}{\log x} +O(1)\\
& \soe \Re \sum_{n\ioe x} \frac{\Lambda (n)}{n^{\sigma +it}\log n} \frac{\log (x/n)}{\log x}- 2 \frac{F(\sigma+it)}{(\log T)/V} +O(1)\\
& \soe -2V -2\frac{V}{\log T}F(\sigma+it)+O(1) \quad \text{\footnotesize (car $t$ est $V$-typique ; on utilise \textit{(i)} page \pageref{t59}).} 
\end{align*}

Pour majorer $F(\sigma+it)$, nous considérons différents domaines de sommation :

$\bullet$ $2\pi n \delta V/\log T \ioe |t-\gamma| \ioe 2\pi (n+1)\delta V /\log T  \quad (0 \ioe n \ioe N=\lfloor (\log T)/4\pi\delta V \rfloor)$,

$\bullet$ le domaine complémentaire, qui est inclus dans $\{ \gamma, \, |t-\gamma|\soe 1/2\}$.

La contribution de la première catégorie de zéros est
\begin{align*}
 & \ioe 2(1+\delta)V \sum_{n=0}^N \frac{(\sigma-1/2)}{(\sigma-1/2)^2+(2\pi\delta nV/\log T)^2}\quad \text{\footnotesize (car $t$ est $V$-typique ; on utilise \textit{(ii)})} \\
& \ll V \Bigl ( \frac{1}{\sigma-1/2} + \frac{1}{\delta V/\log T} \Bigr ) \quad \text{\footnotesize (proposition \ref{t57})} \\
&\ll \delta^{-1}\log T.
\end{align*}

La contribution de la seconde catégorie de zéros est
\begin{align*}
 & \ioe \sum_{|t-\gamma| \soe 1/2}\frac{(\sigma-1/2)}{(\sigma-1/2)^2+(t-\gamma)^2}\\
&\ll \log T \quad \text{\footnotesize (cf. \cite{D2000}, p.98).}
\end{align*}

Cela démontre l'estimation annoncée.\fin

\begin{prop}\label{t61}(HR)
Soit 

$\bullet$ $T$ assez grand ;

$\bullet$  $V$ tel que $(\log \log T)^2 \ioe V \ioe \log T/\log\log T$ ;

$\bullet$ $t$ une ordonnée $V$-typique de taille $T$.

Alors, pour $\demi < \sigma \ioe \sigma_0(=1/2+V/\log T)$, on a 
$$
\log |\zeta (\sigma +it)| \soe \log |\zeta (\sigma_0 +it)| -V\log \frac{(\sigma_0-1/2)}{(\sigma-1/2)}-2(1+\delta)V\log\log V+O(V\delta^{-2}).
$$  
\end{prop}
\dem

On a
\begin{align}\label{t58}
  \log |\zeta (\sigma_0+it)|-\log |\zeta (\sigma+it)| &= \int_{\sigma}^{\sigma_0}\Re \dlz (u+it)du \notag\\
&\ioe \int_{\sigma}^{\sigma_0}F (u+it)du \notag\\
&=\sum_{\gamma}\int_{\sigma}^{\sigma_0}\frac{(u-1/2)}{(u-1/2)^2+(t-\gamma)^2}du \notag\\
&= \demi \sum_{\gamma}\log\frac{(\sigma_0-1/2)^2+(t-\gamma)^2}{(\sigma-1/2)^2+(t-\gamma)^2}.
\end{align}

Pour commencer, observons que
$$
\frac{(\sigma_0-1/2)^2+(t-\gamma)^2}{(\sigma-1/2)^2+(t-\gamma)^2}
$$
est une fonction décroissante de $|t-\gamma|$ si $1/2<\sigma \ioe \sigma_0$. 

Pour majorer la somme \eqref{t58}, nous considérons différents domaines de sommation :

$\bullet$ $|t-\gamma|\ioe \pi V/\bigl ((\log V)\log T\bigr )$ ;

$\bullet$ $(2\pi \delta n+\pi/\log V)V/\log T \ioe |t-\gamma| \ioe \bigl (2\pi \delta (n+1)+\pi/\log V\bigr)V/\log T \quad (0 \ioe n \ioe N=\lfloor (\log T)/4\pi \delta V \rfloor )$ ;

$\bullet$ le domaine complémentaire, qui est inclus dans $\{ \gamma, \, |t-\gamma|\soe 1/2\}$.

Comme on a
$$
\frac{(\sigma_0-1/2)^2+(t-\gamma)^2}{(\sigma-1/2)^2+(t-\gamma)^2} \ioe \frac{(\sigma_0-1/2)^2}{(\sigma-1/2)^2},
$$
et comme $t$ est $V$-typique, la condition \textit{(iii)} (page \pageref{t59}) nous permet de dire que la contribution de la première catégorie de zéros est
$$
\ioe V\log \Bigl (\frac{\sigma_0-1/2}{\sigma-1/2} \Bigr ).
$$

La contribution de la deuxième catégorie de zéros est, d'après \textit{(ii)} (page \pageref{t59})
\begin{align*}
  &\ioe 2(1+\delta)V \cdot \demi \sum_{n=0}^N \log \frac{1+(2\pi \delta n+\pi/\log V)^2}{(2\pi \delta n+\pi/\log V)^2}\\
&= (1+\delta)V\log \bigl (1+(\pi^{-1}\log V)^2 \bigr ) +(1+\delta)V  \sum_{n=1}^N \log \Bigl ( 1+\frac{1}{(2\pi \delta n+\pi/\log V)^2}\Bigr )\\
&\ioe 2(1+\delta)V\log \log V +O(\delta^{-2}V).
\end{align*}

Enfin, la contribution de la troisième catégorie de zéros est
\begin{align*}
  &\ioe\demi \sum_{|t-\gamma| \soe 1/2} \log \Bigl ( 1+\frac{(\sigma_0-1/2)^2}{(t-\gamma)^2} \Bigr )\\
& \ioe\demi \sum_{|t-\gamma| \soe 1/2} \frac{(\sigma_0-1/2)^2}{(t-\gamma)^2} \\
& \ll (V/\log T)^2 \log T\\
& \ll V/\log\log T.
\end{align*}

Ces trois estimations démontrent le résultat annoncé.\fin

\section{Majoration de $|x^z\zeta(z)^{-1}|$}

C'est la proposition suivante qui sera utilisée lors de l'application de la formule de Perron.
\begin{prop}(HR)\label{t43}
Soit 

$\bullet$ $t$ assez grand ;

$\bullet$ $x\soe t$ ;

$\bullet$ $V'$ tel que $(\log \log t)^2 \ioe V' \ioe (\log t/2) /(\log \log t/2)$ ;

$\bullet$ $V \soe V'$.
 
On suppose que $t$ est une ordonnée $V'$-typique (de taille $T'$).  

Alors 
$$
|x^z\zeta(z)^{-1}| \ioe \sqrt{x} \exp \bigl ( V\log (\log x/\log t) +2(1+\delta)V\log \log V +O(V\delta^{-2})\bigr) \quad\quad (V' \ioe (\Re z -1/2)\log x \ioe V, \quad |\Im z|=t).
$$
\end{prop}
\dem

On a 
\begin{align*}
\log |x^z\zeta(z)^{-1}| &= \Re z \log x -\log |\zeta (z)|\\
& \ioe \demi \log x +V -\log |\zeta (z)|,
\end{align*}
car $(\Re z -1/2)\log x \ioe V$. Distinguons maintenant deux cas.

$\bullet$ Si $\Re z -1/2 \ioe V'/\log T'$, la proposition \ref{t40} permet d'affirmer que
\begin{align*}
-\log |\zeta (z)| &\ioe V'\log \frac{V'/\log T'}{\Re z -1/2} +2(1+\delta)V'\log \log V' +O(V'\delta^{-2}) \\
& \ioe V'\log \frac{V'/\log T'}{V'/\log x} +2(1+\delta)V'\log \log V' +O(V'\delta^{-2}) \\
& =  V'\log \frac{\log x}{\log T'} +2(1+\delta)V'\log \log V' +O(V'\delta^{-2}).
\end{align*}

$\bullet$ Si $\Re z -1/2 > V'/\log T'$, on aura grâce à la proposition \ref{t40}
\begin{align*}
-\log |\zeta (z)| &\ioe O(V'\delta^{-1}) \\
& \ioe V'\log \frac{\log x}{\log T'} +2(1+\delta)V'\log \log V' +O(V'\delta^{-2}).
\end{align*}

Ainsi
\begin{align*}
\log |x^z\zeta(z)^{-1}| & \ioe \demi \log x +V -\log |\zeta (z)|\\
& \ioe  \demi \log x + V'\log \frac{\log x}{\log t} +2(1+\delta)V'\log \log V' +O(V'\delta^{-2}) +V +V'\log\frac{\log t}{\log T'}\\
& \ioe  \demi \log x + V\log (\log x/\log t) +2(1+\delta)V\log \log V +O(V\delta^{-2}).\fine
\end{align*}

\section{Préalables à l'estimation de la fréquence des ordonnées atypiques}

Pour la commodité du lecteur, nous rassemblons ici des résultats auxiliaires, essentiellement issus du \S 2 de \cite{S2008}, que nous utiliserons au \S \ref{t29} ci-dessous. Nous renvoyons à \cite{S2008} pour des références et des indications de démonstrations.

\subsection{Encadrement de la fonction indicatrice d'un intervalle}\label{t52}
 
La proposition suivante résume la construction de Selberg de bonnes approximations analytiques de la fonction caractéristique d'un intervalle. On définit la transformation de Fourier de $f \in L^1(\Real)$ par la formule
$$
\hat{f} (x)=\int_{-\infty}^{+\infty}f(t) e^{-2i\pi tx}dt.
$$

\begin{prop}\label{t30}
  Soit $h>0$, $\Delta >0$. Soit $\chi_h=1_{[-h,h]}$ la fonction caractéristique de l'intervalle $[-h,h]$. Il existe des fonctions entières paires $F_+$ et $F_-$ ayant les propriétés suivantes.

\medskip

(i)  $F_-(u) \ioe \chi_h(u) \ioe F_+(u)$ pour tout réel $u$ ;

\medskip

(ii) $\int_{-\infty}^{+\infty} |F_{\pm}(u) - \chi_h(u)|du=1/\Delta$, c'est-à-dire $\hat{F}_{\pm}(0)=2h\pm 1/\Delta$ ;

\medskip

(iii) $\hat{F}_{\pm}$ est réelle et paire, $\hat{F}_{\pm}(x)=0$ pour $|x| \soe \Delta$ et $|x\hat{F}_{\pm}(x)| \ioe 2$ pour tout $x$ ;

\medskip

(iv) Pour $z \in \Com$, $|z| \soe 2h$, on a
$$
|F_{\pm} (z)| \ll \frac{e^{2\pi |\Im z|}}{(\Delta |z|)^2}.
$$
\end{prop}
\dem

Contentons-nous d'une indication pour \textit{(iii)}. On a
\begin{align*}
|x\hat{F}_{\pm}(x)| &\ioe  |x\hat{\chi}_h(x)| +|x|\int_{-\infty}^{+\infty} |F_{\pm}(u) - \chi_h(u)|du\\
& =\Bigl \lvert \frac{\sin 2\pi hx}{\pi} \Bigr \rvert +\frac{|x|}{\Delta}\\
& \ioe \frac{1}{\pi} +1, \quad \text{si $|x| \ioe \Delta$.} 
\end{align*}
D'autre part $|x\hat{F}_{\pm}(x)|=0$ si  $|x| \soe \Delta$.\fin

\subsection{Zéros de la fonction $\zeta$ et nombres premiers : la formule explicite de Guinand-Weil}\label{t51}

Nous donnons une version de la formule explicite de Guinand-Weil, reliant les nombres premiers (représentés par la fonction $\Lambda$ de von Mangoldt) aux zéros non triviaux $\rho=\beta +i\gamma$ de la fonction $\zeta$.

\begin{prop}\label{t31}
  Soit $\eps >0$, $\delta >0$, et $h(s)$ une fonction analytique dans la bande $|\Im s| \ioe \demi +\eps$, y vérifiant $h(s) \ll (1+|s|)^{-1-\delta}$, et réelle sur la droite réelle. Alors
  \begin{align*}
    \sum_{\rho}h\bigl((\rho-i)/2  )\bigr) &= h \Bigl ( \frac{1}{2i} \Bigr )+h \Bigl (- \frac{1}{2i} \Bigr )+\frac{1}{2\pi}\int_{-\infty}^{+\infty}h(u)  \Re \psi \Bigl ( \frac{1}{4}+ \frac{iu}{2} \Bigr ) du\\
&-\frac{\log \pi}{2\pi} \hat{h}(0)- \frac{1}{2\pi} \sum_n \frac{\Lambda (n)}{\sqrt{n}}   \biggl ( \hat{h}\Bigl (\frac{\log n}{2\pi}\Bigr ) +\hat{h}\Bigl (-\frac{\log n}{2\pi}\Bigr )\biggr ).
  \end{align*}
\end{prop}

Signalons tout de suite la formule (2.3) de l'article de Goldston et Gonek \cite{GG2007}. C'est une évaluation de l'intégrale figurant dans la formule explicite, appliquée aux translatées des fonctions $F_{\pm}$.

\begin{prop}\label{t32}
  Soit $t\soe 4$, $\Delta \soe 1$, $0<h \ioe \sqrt{t}$, et $F_{\pm}$ les fonctions de la proposition \ref{t30}. On a
$$
\int_{-\infty}^{+\infty} F_{\pm}(u-t) \Re \psi \Bigl ( \frac{1}{4}+ \frac{iu}{2} \Bigr ) du =(2h \pm 1/\Delta) \log \frac{t}{2} +O(1).
$$
\end{prop}

\subsection{Moments de polynômes de Dirichlet}

Dans la proposition suivante, une inégalité de type {\og grand crible\fg}, nous reprenons en partie la formulation originale de Maier et Montgomery (Lemma 5 de \cite{MM2008}).

\begin{prop}\label{t35}
  Soit $P(s)=\sum_{p\ioe N}a(p)p^{-s}$ un polynôme de Dirichlet (la variable de sommation $p$ est un nombre premier).

Soit $s_1, \dots, s_R \in \Com$, $\alpha \in \Real$ et $T \soe 3$ tels que

$\bullet$ $1 \ioe |\Im (s_i-s_j)| \ioe T$, $i \not =j$ ;

$\bullet$ $\Re s_i \soe \alpha$, $1 \ioe i \ioe R$.

Alors, pour tout entier positif $k$ tel que $N^k \ioe T$, on a
$$
\sum_{r=1}^R |P(s_r)|^{2k} \ll T(\log T)^2k!\Bigl (\sum_{p\ioe N}|a(p)|^2p^{-2\alpha}\Bigr )^k.
$$ 
\end{prop}

\subsection{Un calcul auxiliaire}

Cette proposition est un simple calcul auxiliaire, faisant intervenir quatre paramètres liés par diverses relations.
\begin{prop}\label{t33}
Soit :

$\bullet$ $T$ assez grand ;

$\bullet$ $(\log \log T)^2 \ioe V \ioe \log T/\log\log T$ ;

$\bullet$ $\eta=1/\log V$ ;

$\bullet$ $k=\lfloor  V/(1+\eta)\rfloor$.

Alors
$$
k \bigl ( \log (k\log \log T)-2\log (\eta V) \bigr ) \ioe - V \log (V/\log \log T) +2V \log \log V +V.
$$  
\end{prop}
\dem

On a $k \ioe V$ donc
\begin{align*}
 \log (k\log \log T)-2\log (\eta V) & \ioe  -\log (V/\log \log T) +2 \log \log V\\
& \ioe 0. 
\end{align*}

D'autre part,
\begin{align*}
  k & \soe  V/(1+\eta) -1\\
& \soe  V(1-\eta).
\end{align*}

Par conséquent
\begin{align*}
k \bigl ( \log (k\log \log T)-2\log (\eta V) \bigr ) & \ioe V(1 -\eta )\bigl ( -\log (V/\log \log T) +2 \log \log V \bigr )\\
&\ioe - V \log (V/\log \log T) +2 V \log \log V +\eta V \log V\\
& = -V \log (V/\log \log T) +2 V\log \log V +V.\fine
\end{align*}

\section{Sur le nombre de zéros de la fonction $\zeta$ dans un intervalle de la droite critique}

\subsection{Encadrement paramétrique de l'écart à sa moyenne du nombre de zéros de la fonction $\zeta$ dans un intervalle de la droite critique}\label{t29} 

La proposition suivante donne un encadrement du nombre d'ordonnées de zéros de $\zeta$ dans l'intervalle $]t-h,t+h]$, centré par rapport à sa valeur moyenne $(h/\pi)\log(t/2\pi)$. Cet encadrement est exprimé au moyen d'un paramètre $\Delta$, et met notamment en jeu un polynôme de Dirichlet de longueur $\exp 2\pi \Delta$.
\begin{prop}\label{t34}(HR)
Soit $t \soe 4$, $\Delta \soe 2$ et $0<h \ioe \sqrt{t}$. On a
$$
N(t+h)-N(t-h)- 2h\frac{\log t/2\pi}{2\pi}\ioe  \frac {\log t}{2\pi\Delta}-\frac{1}{\pi}\Re \sum_{p\ioe e^{2\pi\Delta}} \frac{\log p}{p^{\demi+it}} \hat{F}_+\Bigl (\frac{\log p}{2\pi}\Bigr )+O (\log  \Delta ),
$$  
et
$$
N(t+h)-N(t-h)- 2h\frac{\log t/2\pi}{2\pi}\soe - \frac {\log t}{2\pi\Delta}-\frac{1}{\pi}\Re \sum_{p\ioe e^{2\pi\Delta}} \frac{\log p}{p^{\demi+it}} \hat{F}_-\Bigl (\frac{\log p}{2\pi}\Bigr )+O (\log  \Delta ),
$$ 
où $F_+$ et $F_-$ sont les fonctions de la proposition \ref{t30}. De plus, on a pour tout nombre premier $p$ :
$$
\Bigl \lvert \frac{\log p }{\pi}\hat{F}_{\pm}\Bigl (\frac{\log p}{2\pi}\Bigr )\Bigr \rvert\ioe 4.
$$
\end{prop}
\dem

On a
\begin{align*}
  N(t+h)-N(t-h) &= \sum_{\rho}1_{]-h,h]}(\gamma -t)\\
&\ioe \sum_{\rho}\chi_h(\gamma -t)\\
&\ioe  \sum_{\rho}F_+(\gamma -t) \quad \text{\footnotesize (d'après la proposition \ref{t30}, {\it (i)})}\\
&= \sum_{\rho}f(\gamma) \quad \text{\footnotesize (avec $f(u)=F_+(u-t)$)}\\
&=f \Bigl ( \frac{1}{2i} \Bigr )+f \Bigl (- \frac{1}{2i} \Bigr )+\frac{1}{2\pi}\int_{-\infty}^{+\infty}f(u) \Re \psi \Bigl ( \frac{1}{4}+ \frac{iu}{2} \Bigr ) du\\
&-\frac{\log \pi}{2\pi} \hat{f}(0)- \frac{1}{2\pi} \sum_n \frac{\Lambda (n)}{\sqrt{n}}   \biggl ( \hat{f}\Bigl (\frac{\log n}{2\pi}\Bigr ) +\hat{f}\Bigl (-\frac{\log n}{2\pi}\Bigr )\biggr ) \quad\text{\footnotesize (car $f$ vérifie les hypothèses de la proposition \ref{t31})}\\
&=F_+ \Bigl ( \frac{1}{2i}-t \Bigr )+F_+ \Bigl (- \frac{1}{2i}-t \Bigr )+\frac{1}{2\pi}\int_{-\infty}^{+\infty}F_+(u-t)  \Re \psi \Bigl ( \frac{1}{4}+ \frac{iu}{2} \Bigr )du\\
&-\frac{\log \pi}{2\pi} \hat{F}_+(0)- \frac{1}{2\pi} \sum_n \frac{\Lambda (n)}{\sqrt{n}}   \biggl (n^{-it} \hat{F}_+\Bigl (\frac{\log n}{2\pi}\Bigr ) +n^{it}\hat{F}_+\Bigl (-\frac{\log n}{2\pi}\Bigr )\biggr )\quad\text{\footnotesize (car $\hat{f}(x)=e^{-2i\pi tx}\hat{F}_+(x)$)}.
\end{align*}

Examinons successivement les différents termes de cette somme. Nous avons :
\begin{align*}
  F_+ \Bigl ( \frac{1}{2i}-t \Bigr )+F_+ \Bigl (- \frac{1}{2i}-t \Bigr ) &\ll (\Delta t)^{-2} \quad \text{\footnotesize (d'après la proposition \ref{t30} {\it (iv)}, qui s'applique car $t\soe 2h$)}\\
& \ll 1 \quad\text{\footnotesize (car $\Delta \soe 2$ et $t \soe 4$)} ;\\
\frac{1}{2\pi}\int_{-\infty}^{+\infty}F_+(u-t)  \Re \psi \Bigl ( \frac{1}{4}+ \frac{iu}{2} \Bigr )du -\frac{\log \pi}{2\pi} \hat{F}_+(0)&= (2h + 1/\Delta) \frac {\log t/2\pi}{2\pi} +O(1) \quad\text{\footnotesize (propositions \ref{t32} et \ref{t30} {\it (ii)})}\\
\frac{1}{2\pi} \sum_n \frac{\Lambda (n)}{\sqrt{n}}   \biggl (n^{-it} \hat{F}_+\Bigl (\frac{\log n}{2\pi}\Bigr ) +n^{it}\hat{F}_+\Bigl (-\frac{\log n}{2\pi}\Bigr )\biggr ) &=\frac{1}{\pi}\Re \sum_n \frac{\Lambda (n)}{n^{\demi+it}} \hat{F}_+\Bigl (\frac{\log n}{2\pi}\Bigr ) \quad\text{\footnotesize (car $\hat{F}_+$ est réelle et paire)}\\
&=\frac{1}{\pi}\Re \sum_{n\ioe e^{2\pi\Delta}} \frac{\Lambda (n)}{n^{\demi+it}} \hat{F}_+\Bigl (\frac{\log n}{2\pi}\Bigr )\quad\text{\footnotesize (car $\hat{F}_+(x)=0$ pour $x \soe \Delta$)}\\
&= \frac{1}{\pi}\Re \sum_{p\ioe e^{2\pi\Delta}} \frac{\log p}{p^{\demi+it}} \hat{F}_+\Bigl (\frac{\log p}{2\pi}\Bigr )+\frac{1}{\pi}\Re \sum_{p\ioe e^{\pi\Delta}} \frac{\log p}{p^{1+2it}} \hat{F}_+\Bigl (\frac{\log p}{\pi}\Bigr )+O(1)\\
& \quad\text{\scriptsize (car $\Lambda (n)\hat{F}_+(\log n /2\pi)\ll 1$ ; la contribution des $n=p^k$ avec $k \soe 3$ est donc $O(1)$)}\\
&= \frac{1}{\pi}\Re \sum_{p\ioe e^{2\pi\Delta}} \frac{\log p}{p^{\demi+it}} \hat{F}_+\Bigl (\frac{\log p}{2\pi}\Bigr )+O (\log  \Delta ).
 \end{align*}

Finalement, nous avons établi la majoration
$$
N(t+h)-N(t-h)- 2h\frac{\log t/2\pi}{2\pi}\ioe  \frac {\log t}{2\pi\Delta}-\frac{1}{\pi}\Re \sum_{p\ioe e^{2\pi\Delta}} \frac{\log p}{p^{\demi+it}} \hat{F}_+\Bigl (\frac{\log p}{2\pi}\Bigr )+O (\log  \Delta ).
$$
La démonstration de la minoration est analogue.\fin

\medskip

\subsection{Majoration de l'écart à sa moyenne du nombre de zéros de la fonction $\zeta$ dans un intervalle de la droite critique}

Lorsqu'on majore trivialement les polynômes de Dirichlet qui interviennent dans la proposition \ref{t34}, on obtient le résultat suivant, dû à Goldston et Gonek (cf. \cite{GG2007}). Notre énoncé est légèrement plus précis que celui de \cite{GG2007}.
\begin{prop}\label{t42}
 Soit $t$ assez grand et $0<h \ioe \sqrt{t}$. On a
$$
|N(t+h)-N(t-h)- (h/\pi)\log (t/2\pi)|\ioe (\log t)/2\log \log t + \bigl (1/2 +o(1)\bigr ) \log t\log \log \log t/(\log \log t)^2.
$$ 
\end{prop}
\dem

On a
\begin{align*}
  \Bigl \lvert \frac{1}{\pi} \sum_{p\ioe e^{2\pi\Delta}} \frac{\log p}{p^{\demi+it}} \hat{F}_{\pm}\Bigl (\frac{\log p}{2\pi}\Bigr )\Bigr \rvert& \ll \sum_{p\ioe e^{2\pi\Delta}} \frac{1}{\sqrt{p}}\\
&\ll \frac{ e^{\pi\Delta}}{\Delta}.
\end{align*}
On choisit $\Delta = \frac{1}{\pi}\log(\log t/\log \log t)$ et on vérifie alors que 
\begin{equation*}
 \frac {\log t}{2\pi\Delta}+O(e^{\pi\Delta}/\Delta) +O (\log  \Delta )= (\log t)/2\log \log t +  \bigl (1/2 +o(1)\bigr ) \log t\log \log \log t/(\log \log t)^2.\fine 
\end{equation*}

\subsection{Fréquence des déviations du nombre de zéros de la fonction $\zeta$ dans un intervalle de la droite critique}
Nous donnons maintenant une majoration du nombre de points $t$ {\og bien espacés\fg} de l'intervalle $[T,2T]$ pour lesquels le nombre d'ordonnées de zéros de $\zeta$ dans l'intervalle $]t-h,t+h]$ dépasse sa moyenne $(h/\pi)\log(t/2\pi)$ d'une quantité $V$. On a une majoration analogue pour la fréquence des déviations dans l'autre direction, mais nous n'en aurons pas l'usage.

\begin{prop}\label{t37}
Soit :

$\bullet$ $T$ assez grand ;

$\bullet$ $0<h \ioe \sqrt{T}$ ;

$\bullet$ $(\log \log T)^2 \ioe V \ioe \log T/\log\log T$ ;

$\bullet$ $T\ioe t_1<t_2< \dots <t_R\ioe 2T$ tels que $t_{r+1}-t_r \soe 1$, $1 \ioe r <R$.

On suppose que
$$
N(t_r+h)-N(t_r-h)- h\frac{\log t_r/2\pi}{\pi} \soe V, \quad 1 \ioe r \ioe R.
$$
Alors
$$
R \ll T\exp \bigl (- V \log (V/\log \log T) +2V \log \log V +O(V)\bigr ).
$$
\end{prop}
\dem

La majoration de $N(t_r+h)-N(t_r-h)- (h/\pi)\log (t_r/2\pi)$ fournie par la proposition \ref{t34} montre que pour $\Delta \soe 2$, on a
$$
 \Bigl \lvert \sum_{p\ioe e^{2\pi\Delta}} \frac{a(p)}{p^{\demi+it_r}} \Bigr \rvert \soe V-\frac{\log 2T}{2\pi \Delta} +O(\log \Delta), \quad 1\ioe r \ioe R,
$$
où $a(p) =a(p, \Delta,h)=\pi^{-1}\log p\, \hat{F}_+\bigl ((\log p)/2\pi \bigr )$ vérifie $|a(p)| \ioe 4$, d'après la proposition \ref{t30}, \textit{(iii)}.

On choisit 
$$
\Delta = \frac{(1+\eta)\log T}{2\pi V} \quad \text{avec $\eta =1/\log V$},
$$
de sorte que
\begin{align*}
 \exp 2\pi \Delta &=T^{(1+\eta)/V},\\
\log \Delta &\ll \log \log T \ll \sqrt{V},
\end{align*}
et
\begin{align*}
V-\frac{\log 2T}{2\pi \Delta} +O(\log \Delta) &= \eta \frac{V}{1+\eta}+O(\sqrt{V})\\
& \soe \demi \eta V  
\end{align*}

Par conséquent,
$$
\Bigl \lvert \sum_{p\ioe T^{(1+\eta)/V}} \frac{a(p)}{p^{\demi+it_r}} \Bigr \rvert \soe \demi \eta V, \quad 1\ioe r \ioe R.
$$
En élevant cette inégalité à la puissance $2k$ et en sommant pour $r=1, \dots, R$, on obtient
\begin{align*}
  R(\eta V/2)^{2k} & \ioe \sum_{r=1}^R \Bigl \lvert \sum_{p\ioe T^{(1+\eta)/V}} \frac{a(p)}{p^{\demi+it_r}} \Bigr \rvert^{2k}\\
& \ll T(\log T)^2k!\Bigl (\sum_{p\ioe T^{(1+\eta)/V}}|a(p)|^2p^{-1}\Bigr )^k,
\end{align*}
d'après la proposition \ref{t35}, pourvu que $T^{k(1+\eta)/V}\ioe T$.

La dernière quantité est 
$$
\ll T(\log T)^2 (Ck\log \log T)^k,
$$
où $C$ est une constante absolue. Ainsi
$$
R \ll T(\log T)^2 (4C)^k\bigl ( (k\log \log T)/\eta^2V^2\bigr )^k.
$$
On choisit $k=\lfloor  V/(1+\eta) \rfloor$. La proposition \ref{t33} s'applique :
$$
( (k\log \log T)/\eta^2V^2\bigr )^k \ioe  \exp \bigl (-  V \log (V/\log \log T) +2V \log \log V+V \bigr ).
$$
Enfin,
\begin{align*}
(\log T)^2 (4C)^k  &= \exp (k\log 4C +2 \log \log T)\\
&=\exp O(V),  
\end{align*}
d'où le résultat.\fin

\medskip

\section{Valeurs de $V$ pour lesquelles toutes les ordonnées sont $V$-typiques}

Nous donnons une variante un peu plus précise de la première assertion de la Proposition 4 de \cite{S2008}.
\begin{prop}\label{t41}
  Soit $T$ assez grand, et $V$ tel que
$$
\demi + \log \log \log T/\log \log T \ioe V \log \log T/ \log T \ioe 1.
$$
Alors toute ordonnée $t \in [T,2T]$ est $V$-typique.
\end{prop}
\dem

Il faut vérifier les critères \textit{(i), (ii), (iii)} de la définition d'une ordonnée $V$-typique.

Pour \textit{(i)}, on a
$$
f(u)=\sum_{n\ioe u}\frac{\Lambda (n)}{\sqrt{n}\log n} \ll \frac{\sqrt{u}}{\log u}, \quad u \soe 2,
$$
donc
\begin{align*}
  \sum_{n\ioe x}\frac{\Lambda (n)}{\sqrt{n}\log n}\log (x/n)  &= \int_1^x f(u)\frac{du}{u}\\
& \ll \frac{\sqrt{x}}{\log x}, \quad x \soe 2.
\end{align*}

Or $T^{1/V} \ioe (\log T)^2$ car $V \soe \demi \log T/\log \log T $. Par conséquent, pour $\sigma \soe 1/2$, $t \in \Real$, et $x=T^{1/V}$, on a
\begin{align*}
\Bigl \lvert \sum_{n \ioe x} \frac{\Lambda (n)}{n^{\sigma +it}\log n}\frac{\log (x/n)}{\log x} \Bigr \rvert & \ioe  \sum_{n\ioe x}\frac{\Lambda (n)}{\sqrt{n}\log n}\frac{\log (x/n)}{\log x}\\
& \ll  \frac{\sqrt{x}}{(\log x)^2}\\
& \ll  \frac{\log T}{(\log \log T)^2}\\
& = o(V). 
\end{align*}
 
Pour \textit{(ii)} on a, avec $t' \in[t-1,t+1]$ et $h=\pi\delta V/\log T$ :
\begin{align*}
1+N(t'+h)-N(t'-h) &\ioe  (h/\pi)\log (t'/2\pi) +\demi \log t'/\log \log t' + \bigl (1/2 +o(1)\bigr ) \log t'\log \log \log t'/(\log \log t')^2\\
 & \quad \text{\footnotesize (proposition \ref{t42})}\\ 
& \ioe  (h/\pi)\log T + \demi \log T/\log \log T + \log T\log \log \log T/(\log \log T)^2\\
& \ioe (1+\delta)V.
\end{align*}

Pour \textit{(iii)} on a, avec $t' \in[t-1,t+1]$ et $h=\pi V/\bigl ((\log V)\log T\bigr )$ :
\begin{align*}
1+N(t'+h)-N(t'-h) &\ioe  (h/\pi)\log (t'/2\pi) +\demi \log t'/\log \log t' + \bigl (1/2 +o(1)\bigr ) \log t'\log \log \log t'/(\log \log t')^2\\
 & \ioe  \frac{V}{\log V} + \demi \log T/\log \log T +  \bigl (1/2 +o(1)\bigr )\log T\log \log \log T/(\log \log T)^2\\
& \ioe \demi \log T/\log \log T + \log T\log \log \log T/(\log \log T)^2\\
& \ioe V.\fine
\end{align*}

\subsection{Minoration du logarithme du module de la fonction $\zeta$}

\begin{prop}\label{t62}
  Pour $|t|$ assez grand et $1/2<\sigma \ioe 2$, on a
$$
\log |\zeta (\sigma +it)| \soe -\frac{\log |t|}{\log \log |t|}\log \frac{1}{\sigma-1/2}-3\frac{\log |t|\log \log \log |t|}{\log \log |t|}.
$$
\end{prop}
\dem

On applique les propositions \ref{t40} et \ref{t41} en prenant
$$
V= \frac{\log |t|}{\log \log |t|}, \quad \delta=\demi.
$$

L'ordonnée $|t|$ est $V$-typique de taille $|t|$ et
\begin{align*}
V\log \Bigl( \frac{V/\log |t|}{\sigma-1/2} \Bigr )+2(1+\delta)V\log\log V +O(V\delta^{-2}) &=\frac{\log |t|}{\log \log |t|}\log \frac{1}{\sigma-1/2}+2\frac{\log |t|\log \log \log |t|}{\log \log |t|} +O(\log |t|/\log \log |t|)\\
& \ioe \frac{\log |t|}{\log \log |t|}\log \frac{1}{\sigma-1/2}+3\frac{\log |t|\log \log \log |t|}{\log \log |t|}.\fine    
\end{align*}

\section{Majoration du nombre d'ordonnées atypiques bien espacées dans un intervalle}\label{t39}

Nous majorons maintenant le nombre d'ordonnées $V$-atypiques de taille $T$, bien espacées. 

\begin{prop}(RH)\label{t64}
Soit

$\bullet$ $T$ assez grand ;

$\bullet$ $2(\log \log T)^2 \ioe V \ioe \log T/\log\log T$ ;

$\bullet$ $T\ioe t_1<t_2< \dots <t_R\ioe 2T$ des ordonnées $V$-atypiques telles que $t_{r+1}-t_r \soe 1$, $1 \ioe r <R$.

Alors
$$
R \ll T\exp \bigl (- V \log (V/\log \log T) +2 V\log \log V +O(V)\bigr ).
$$  
\end{prop}
\dem

Observons d'abord que le majorant annoncé est une fonction décroissante de $V$ et que sa valeur en $V=\log T/\log \log T$ est 
$$
\soe \exp (3 \log T \log \log \log T/\log \log T),
$$
quantité qui tend vers l'infini avec $T$.

Il suffit de démontrer séparément la majoration pour le nombre $R_1$ (resp. $R_2$, resp. $R_3$) (mais nous le noterons encore $R$) de points (que nous noterons encore $t_1, \dots, t_R$) bien espacés dans $[T,2T]$ infirmant la condition \textit{(i)} (resp. \textit{(ii)}, resp. \textit{(iii)}) de la définition page \pageref{t59}.

Commençons par la condition \textit{(i)}. Si elle est en défaut pour chaque $t_r$, il existe des $\sigma_r \soe 1/2$ tels que
$$
\Bigl \lvert \sum_{n \ioe x} \frac{\Lambda (n)}{n^{\sigma_r +it_r}\log n}\frac{\log (x/n)}{\log x} \Bigr \rvert > 2V, \quad 1 \ioe r \ioe R \quad (x=T^{1/V}).
$$

La contribution des $n=p^{\alpha}$ avec $\alpha \soe 2$ est
\begin{align*}
  &\ll \log \log x\\
& \ll \log \log T\\
& \ll \sqrt{V}.
\end{align*}

Par conséquent,
$$
\Bigl \lvert \sum_{p \ioe x} \frac{1}{p^{\sigma_r +it_r}}\frac{\log (x/p)}{\log x} \Bigr \rvert \soe V, \quad 1 \ioe r \ioe R \quad (x=T^{1/V}).
$$

d'où
\begin{align*}
  RV^{2k} & \ioe \sum_{r=1}^R \Bigl \lvert \sum_{p \ioe x} \frac{1}{p^{\sigma_r +it_r}}\frac{\log (x/p)}{\log x} \Bigr \rvert^{2k}\\
& \ll T(\log T)^2k!\Bigl ( \sum_{p \ioe x} \frac{1}{p}\frac{\log ^2(x/p)}{\log^2 x} \Bigr )^k,
\end{align*}
pourvu que $x^k \ioe T$, c'est-à-dire $k \ioe  V$.

On a
\begin{align*}
\sum_{p \ioe x} \frac{1}{p}\frac{\log ^2(x/p)}{\log^2 x} & \ll \log \log x\\
&\ioe \log \log T,
\end{align*}
et donc
$$
R \ll T(\log T)^2  ( (Ck\log \log T)/V^2\bigr )^k,
$$
où $C$ est une constante absolue. Le choix $k=\lfloor  V \rfloor$ et un calcul plus simple que celui de la proposition \ref{t33} conduisent à l'estimation
$$
R \ll T\exp \bigl ( - V  \log (V/\log \log T) + O(V)\bigr ).
$$

Passons à la condition \textit{(ii)}. Si $t_r$ ne la vérifie pas, il existe $t'_r$ tel que $|t_r-t'_r|\ioe 1$ et
$$
N(t'_r+\pi\delta V/\log T)-N(t'_r-\pi\delta V/\log T)\soe (1+\delta)V-1,
$$
d'où
$$
N(t'_r+\pi\delta V/\log T)-N(t'_r-\pi\delta V/\log T)-(\delta V/\log T)\log (t'_r/2\pi) \soe V +O(1).
$$

Choisissons une ordonnée $t_r$ sur trois, de sorte que les $t'_r$ correspondants soient bien espacés, et gardons uniquement les $t'_r$ de l'intervalle $[T,2T]$ (nous en laissons alors au plus deux de côté). Nous pouvons appliquer la proposition \ref{t37} :
$$
\lfloor (R+2)/3 \rfloor -2 \ll T\exp \bigl (- V \log (V/\log \log T) +2V \log \log V +O(V)\bigr ),
$$
et on a la même estimation pour $R$.

Enfin, si $t_r$ infirme la condition \textit{(iii)}, il existe $t'_r$ tel que $|t_r-t'_r|\ioe 1$ et
$$
N\Bigl (t'_r+\pi V/\bigl ( (\log V)\log T \bigr )\Bigr )-N\Bigl (t'_r-\pi V/\bigl ( (\log V)\log T \bigr )\Bigr )\soe V-1,
$$
d'où
$$
N\Bigl (t'_r+\pi V/\bigl ( (\log V)\log T \bigr )\Bigr )-N\Bigl (t'_r-\pi V/\bigl ( (\log V)\log T \bigr )\Bigr )-V\log (t'_r/2\pi)/\bigl ( (\log V)\log T \bigr )\soe V+O(V/\log V).
$$

En procédant comme pour la condition \textit{(ii)}, nous appliquons de nouveau la proposition \ref{t37}. Il reste à observer que, pour $V'=V+O(V/\log V)$, on a
$$
- V' \log (V'/\log \log T) +2 V' \log \log V' +O(V') \ioe - V \log (V/\log \log T) +2V \log \log V +O(V).
$$
C'est un calcul facile, qui termine la démonstration.\fin

\section{Démonstration du théorème}

\subsection{Application de la formule de Perron}

Rappelons dans l'énoncé suivant la forme classique de la formule de Perron que nous allons utiliser (pour une démonstration, voir \cite{T2008}, II.3).
\begin{prop}\label{t38}
  Soit
$$
F(z)= \sum_{n=1}^{+\infty} \frac{a_n}{n^z},
$$
une s\'erie de Dirichlet, absolument convergente pour $\Re z > 1$. On suppose que $|a_n| \ioe 1$ pour tout $n$. Alors, pour $N \soe T \soe 3$, $N$ entier, on a
$$
\sum_{n \ioe N}a_n = \frac{1}{2 \pi i}\int_{c-iT}^{c+iT} F(z) \frac{N^z}{z}dz + O(N \log T/T),
$$
o\`u $c = 1 +1/ \log N$. La constante implicite dans le O est absolue.
\end{prop}

On a donc 
$$
M(N)=A_N+O(\log N) \quad (N \soe 3),
$$
où
$$
A_N=\frac{1}{2 \pi i}\int_{1+1/ \log N-i2^{\lfloor \log N/\log 2\rfloor}}^{1+1/ \log N+i2^{\lfloor \log N/\log 2\rfloor}} \zeta(z)^{-1} \frac{N^z}{z}dz.
$$

\subsection{Déformation du chemin d'intégration}

Pour majorer $|A_N|$, nous allons remplacer le segment d'intégration $[1+1/\log N-i2^{\lfloor \log N/\log 2\rfloor},1+1/\log N +i2^{\lfloor \log N/\log 2\rfloor}]$ par une variante $\Scal_N$ du chemin défini par Soundararajan dans \cite{S2008}. Nous commençons par une description de $\Scal_N$. On suppose $N$ assez grand. Nous posons
$$
\kappa=\lfloor  (\log N)^{1/2}(\log \log N)^c \rfloor, \quad K =\lfloor \log N/\log 2\rfloor,
$$
où $c$ est une constante positive, fixée ultérieurement. Nous posons également $T_k=2^{k}$ pour $\kappa \ioe k \ioe K$.

\medskip

Le chemin $\Scal_N$ est symétrique par rapport à l'axe réel, et constitué de segments verticaux et horizontaux. Nous décrivons seulement la partie de $\Scal_N$ située dans le demi-plan $\Im z \soe 0$.

$\bullet$ Il y a d'abord un segment vertical $[1/2+1/\log N,1/2+1/\log N +iT_{\kappa}]$.

$\bullet$ Pour chaque $k$ tel que $\kappa \ioe k < K$, on considère les entiers $n$ de l'intervalle $[T_k,2T_k[$. On définit alors $V_n$ comme le plus petit entier de l'intervalle $[(\log \log T_k)^2,\log T_k/\log \log T_k]$ tel que tous les points de $[n,n+1]$ soient $V_n$-typiques de taille $T_k$. L'existence de $V_n$ est garantie par la proposition \ref{t41}. On a même 
$$
V_n \ioe \demi  \log n/\log \log n +\log n (\log \log \log n )/(\log \log n)^2 +1.
$$ 
On inclut alors dans $\Scal_N$ le segment vertical $[1/2+V_n/\log N +in,1/2+V_n/\log N +i(n+1)]$

Il y a enfin des segments horizontaux reliant tous ces segments verticaux :

$\bullet$ le segment $[1/2+1/\log N +iT_{\kappa},1/2+V_{N_0}/\log N +iT_{\kappa}]$ ;

$\bullet$ les segments $[1/2+V_n/\log N +i(n+1),1/2+V_{n+1}/\log N +i(n+1)]$, $T_{\kappa} \ioe n \ioe T_{K}-2$ ;

$\bullet$ le segment $[1/2+V_{N-1}/\log N +iT_K,1+1/\log N +iT_K]$.

\medskip

D'après le théorème de Cauchy, on a sous l'hypothèse de Riemann
$$
A_N=\intc_{\Scal_N}\zeta(z)^{-1}\frac{N^{z}}{z}dz.
$$

\subsection{Majorations des contributions à $A_N$ des différents segments de $\Scal_N$}

\begin{prop}\label{t63}
  On a
$$
 N^{-1/2}A_N \ll_{\delta} \exp \bigl ((\log N)^{1/2} (\log \log N)^{c+\delta}\bigr )  +B_N,
$$
où
$$
B_N=\sum_{n=T_{\kappa}}^{T_K-1}\frac{1}{n} \exp \bigl ( V_n\log (\log N/\log n) +2(1+2\delta)V_n\log\log V_n \bigr ).
$$
\end{prop}
\dem

Pour commencer, nous avons
\begin{align*}
  \biggl \lvert \intc_{\Scal_N, |\Im z|\ioe T_{\kappa}}\zeta(z)^{-1}\frac{N^{z}}{z}dz \biggr \rvert & \ll N^{1/2}\int_{-T_{\kappa}}^{T_{\kappa}}|\zeta(1/2+1/\log N+i\tau)|^{-1} \frac{d\tau}{\sqrt{1/4+\tau^2}}\\
&\ll  N^{1/2}\log T_{\kappa} \max_{|\tau|\ioe T_{\kappa}}|\zeta(1/2+1/\log N+i\tau)|^{-1}\\
&\ioe N^{1/2}(\log T_{\kappa}) \exp \bigl ((\log T_{\kappa}/\log\log T_{\kappa})\log\log N +3(\log T_{\kappa}/\log\log T_{\kappa})\log\log \log T_{\kappa} \bigr )  \\
&\quad \text{\footnotesize (d'après la proposition \ref{t62})}\\
& \ioe N^{1/2}T_{\kappa}^3 \quad \text{\footnotesize (car $\log T_{\kappa} \soe (\log N)^{1/2}$)}.
\end{align*}

\medskip

Pour les autres segments, nous examinons seulement la partie de $\Scal_N$ située dans le demi-plan $\Im z >0$. On aura les mêmes estimations pour la partie située dans le demi-plan $\Im z <0$. 

La contribution du segment horizontal $[1/2+1/\log N \pm iT_{\kappa},1/2+V_{T_{\kappa}}/\log N \pm iT_{\kappa}]$  est
\begin{align*}
 & \ll N^{1/2} (\exp V_{T_{\kappa}}) T_{\kappa}^{-1} \exp \bigl ((\log T_{\kappa}/\log\log T_{\kappa})\log\log N +3(\log T_{\kappa}/\log\log T_{\kappa})\log\log \log T_{\kappa} \bigr )  \\
&\ll N^{1/2}T_{\kappa}^3.
\end{align*}

\medskip

Pour chaque $n$, $T_{\kappa} \ioe n \ioe T_K-1$, la contribution du segment vertical $[1/2+V_n/\log N +in,1/2+V_n/\log N +i(n+1)]$ est
$$
\ll \frac{1}{n}N^{1/2} \exp \bigl ( V_n\log (\log N/\log n) +2(1+\delta)V_n\log\log V_n +DV_n\delta^{-2}\bigr )
$$
(où $D$ est une constante positive absolue), d'après la proposition \ref{t43} (avec $V=V'=V_n$, $n \ioe t \ioe n+1$, $T'=n$, $x=N$).

Pour chaque $n$, $T_{\kappa} \ioe n \ioe T_K-2$, la contribution du segment horizontal $[1/2+V_n/\log N +i(n+1),1/2+V_{n+1}/\log N +i(n+1)]$ est
$$
\ll \frac{1}{n}N^{1/2} \exp \Bigl ( V\log \bigl(\log N/\log (n+1)\bigr ) +2(1+\delta)V\log\log V +DV\delta^{-2}\Bigr ),
$$
où l'on a posé $V=\max(V_n,V_{n+1})$, toujours d'après la proposition \ref{t43}, qui s'applique car $n+1$ est à la fois $V_n$-typique et $V_{n+1}$-typique. Observons que la borne obtenue est
\begin{multline*}
\ll \frac{1}{n}N^{1/2} \exp \bigl ( V_n\log (\log N/\log n) +2(1+\delta)V_n\log\log V_n +DV_n\delta^{-2}\bigr )+\\
\frac{1}{n+1}N^{1/2} \exp \Bigl ( V_{n+1}\log \bigl (\log N/\log (n+1) \bigr ) +2(1+\delta)V_{n+1}\log\log V_{n+1} +DV_{n+1}\delta^{-2}\Bigr ).  
\end{multline*}

Enfin, la contribution du  segment horizontal $[1/2+V_{T_K-1}/\log N +iT_K,1+1/\log N +iT_K]$ est
\begin{align*}
  &\ll \frac{1}{T_K}N\exp\bigl (O(\delta^{-1}V_{T_K-1})\bigr ) \quad \text{\footnotesize (d'après la proposition \ref{t40})}\\
& \ll_{\delta} N^{1/2} \quad \text{\footnotesize (car $V_{T_K-1} \ioe \log N/\log\log N$).}\\ 
\end{align*}

En notant en outre que 
$$
\exp(DV\delta^{-2}) \ll_{\delta} \exp(\delta V\log \log V),
$$
et
\begin{align*}
T_{\kappa}^3 &\ioe \exp \bigl (3(\log 2)(\log N)^{1/2} (\log \log N)^{c}\bigr )\\
&   \ll_{\delta} \exp \bigl ((\log N)^{1/2} (\log \log N)^{c+\delta}\bigr ),  
\end{align*}
cela termine la démonstration de la proposition.\fin

\subsection{\'Etude de la somme $B_N$}

Afin de majorer $B_N$, nous allons faire appel au calcul auxiliaire suivant.

\begin{prop}\label{t66}
  Soit $A$ et $C$ des paramètres positifs tels que $A \soe 4C^4+1$. On a alors
$$
AV-V\log V+CV\log\log V \ioe e^AA^C \quad (V>e^C).
$$
\end{prop}
\dem

Posons
$$
f(V)=AV-V\log V+CV\log\log V \quad (V>1).
$$
On a
\begin{align*}
  f'(V) &= A-\log V+C\log \log V-1 +\frac{C}{\log V},\\
f''(V) &= -\frac{1}{V} +\frac{C}{V\log V}-\frac{C}{V(\log V)^2}.
\end{align*}

En particulier, on a $f''(V) <0$ si $V>e^C$. De plus,
\begin{align*}
  f'(e^C) &= A-C+C\log C-1+1\\
&= A-C +C \log C\\
& \soe 4C^4+1-C +C \log C\\
& >0, 
\end{align*}
et $f'(\infty)=-\infty$.

Il existe donc un unique $V_0 >e^C$ tel que

$\bullet$ $f'(V)>0$ pour $e^C \ioe V < V_0$ ;

$\bullet$ $f'(V)<0$ pour $V > V_0$.

On a $f'(V_0)=0$, c'est-à-dire
$$
A-\log V_0+C\log \log V_0=1 -\frac{C}{\log V_0}.
$$

D'autre part,
\begin{align*}
  \max_{V\soe e^C}f(V) &= f(V_0)\\
&= V_0(A-\log V_0 +C \log \log V_0)\\
&=V_0(1-C/\log V_0)\\
& \ioe V_0.
\end{align*}

Posons maintenant $V_1=e^AA^C$ $(>e^C)$. On a
\begin{align*}
  f'(V_1) &=  A-\log V_1+C\log \log V_1-1 +\frac{C}{\log V_1}\\
&= A- (A +C\log A)+C\log (A+C\log A)-1 +\frac{C}{A+C\log A }\\
&\ioe C\log \Bigl ( 1 + \frac{C\log A}{A} \Bigr )-1+\frac{C}{A} \\
& \ioe C\log \Bigl ( 1 + \frac{C}{A^{1/2}} \Bigr )-1 +\frac{C}{A} \quad\text{\footnotesize ( car $\log A \ioe A^{1/2}$)}\\
& \ioe \frac{C^2}{A^{1/2}}-1 +\frac{C}{A}\\
& \ioe 0.
\end{align*}

On a donc $V_0\ioe V_1$, d'où le résultat annoncé.\fin

\begin{prop}\label{t67}
  On a
$$
B_N \ll_{\delta} \exp \bigl ( (\log N)^{1/2} (\log \log N)^{5-c+6\delta}\bigr ).
$$
\end{prop}
\dem

Nous allons utiliser un découpage dyadique en considérant les sommes
$$
B_N(T_k)=\sum_{T_k \ioe n < 2T_k}\frac{1}{n} \exp \bigl ( V_n\log (\log N/\log n) +2(1+2\delta)V_n\log\log V_n \bigr ) \quad (\kappa \ioe k < K).
$$

On aura en effet
\begin{align*}
B_N &\ioe K\max_{\kappa \ioe k < K}B_N(T_k)\\
&\ll \log N \max_{\kappa \ioe k < K}B_N(T_k).
\end{align*}

Posons $T_k=T$. Nous réarrangeons $B_N(T)$ suivant les valeurs de $V_n$.
\begin{align*}
  B_N(T) &= \sum_{\substack{(\log\log T)^2 \ioe V\\ V \ioe (\log T)/(\log \log T)}}\sum_{\substack{T \ioe n < 2T\\ V_n=V}}\frac{1}{n} \exp \bigl ( V\log (\log N/\log n) +2(1+2\delta)V\log\log V \bigr )\\
& \ioe \frac{1}{T}\sum_{\substack{(\log\log T)^2 \ioe V\\ V \ioe (\log T)/(\log \log T)}}\exp \bigl ( V\log (\log N/\log T) +2(1+2\delta)V\log\log V \bigr ) \card \{n, \, T \ioe n < 2T, \, V_n=V\}.
\end{align*}

Nous majorons d'abord la contribution des $V \ioe 2(\log \log T)^2+1$, en utilisant la majoration triviale
$$
\card \{n, \, T \ioe n < 2T, \, V_n=V\} \ioe T.
$$
Cette contribution est
$$
\ioe \exp \Bigl ( O \bigl ( (\log \log N)^3\bigr ) \Bigr ).
$$

\medskip

Considérons maintenant  la contribution des $V > 2(\log \log T)^2+1$. Si $T \ioe n < 2T$ et $V_n=V$, la minimalité de $V_n$ entraîne l'existence dans l'intervalle $[n,n+1]$ de $t_n$, ordonnée $(V_n-1)$-atypique de taille $T$. En évitant de prendre des nombres consécutifs, on partitionne alors l'ensemble
$$
\{n, \, T \ioe n < 2T, \, V_n=V\}
$$
en (au plus) deux sous-ensembles, chacun étant de même cardinal qu'un ensemble d'ordonnées $(V_n-1)$-atypiques de taille $T$, bien espacées dans $[T,2T]$. La proposition \ref{t64} donne alors 
\begin{align*}
 \card \{n, \, T \ioe n < 2T, \, V_n=V\} &\ll T\exp \Bigl (- (V-1) \log \bigl( (V-1)/\log \log T \bigr) +2 (V-1)\log \log (V-1) +O(V)\Bigr )\\
& \ll T\exp \bigl (- V \log ( V/\log \log T ) +2 V\log \log V +O(V)\bigr ). 
\end{align*}

Par conséquent,
\begin{multline}
  \label{t65}
B_N(T) \ll_{\delta} \exp \Bigl ( O \bigl ( (\log \log N)^3\bigr ) \Bigr ) +\\ \sum_{\substack{2(\log\log T)^2 +1 \ioe V\\ V \ioe (\log T)/(\log \log T)}}\exp \Bigl ( V\log \bigl (\log N(\log \log T)/\log T\bigr )-V\log V +(4+5\delta)V\log\log V \Bigr ).  
\end{multline}

\medskip

Pour majorer la somme intervenant dans \eqref{t65}, on utilise la proposition \ref{t66} avec
$$
A=\log \bigl (\log N(\log \log T)/\log T\bigr )\quad \text{et} \quad C=4+5\delta
$$
(on a bien $A \soe 4C^4+1$ et $V>e^C$). 

Par conséquent,
\begin{align*}
\sum_{\substack{2(\log\log T)^2 +1 \ioe V\\ V \ioe (\log T)/(\log \log T)}}&\exp \Bigl ( V\log \bigl (\log N(\log \log T)/\log T\bigr )-V\log V +(4+5\delta)V\log\log V \Bigr )\\
 &\ioe \frac{\log T}{\log \log T}\exp \biggl (\log N\frac{\log \log T}{\log T} \Bigl(\log\bigl ( \log N\frac{\log \log T}{\log T}\bigr )\Bigr)^{4+5\delta}\biggr ).
\end{align*}

Comme
\begin{align*}
 \frac{\log \log T}{\log T} &=\frac{\log \log T_k}{\log T_k}\\
& \ioe \frac{\log \log T_{\kappa}}{\log T_{\kappa}}\\
& \ioe \frac{\log \log N}{(\log N)^{1/2}(\log \log N)^c}, 
\end{align*}
on a finalement
\begin{equation*}
B_N \ll_{\delta} \exp \bigl ( (\log N)^{1/2} (\log \log N)^{5-c+6\delta}\bigr ).\fine  
\end{equation*}

\subsection{Conclusion}

En réunissant les résultats des propositions \ref{t38}, \ref{t63} et \ref{t67}, on obtient 
$$
  N^{-1/2}M(N)  \ll_{\delta} \exp \bigl ((\log N)^{1/2} (\log \log N)^{c+\delta}\bigr ) + \exp \bigl ( (\log N)^{1/2} (\log \log N)^{5-c+6\delta}\bigr ). 
$$

On choisit alors $c=5/2$ et $\delta$ arbitrairement petit pour démontrer le théorème.

\medskip

\noindent Michel BALAZARD\\
C.N.R.S., Institut de Math\'ematiques de Luminy\\
Case 907\\
13288 Marseille Cedex 09\\
FRANCE

\medskip

\noindent Adresse \'electronique : \texttt{balazard@iml.univ-mrs.fr}

\medskip

\noindent Anne de ROTON\\
Institut Elie Cartan, Université Henri Poincaré Nancy 1\\
BP 239\\
54506 Vandoeuvre-lès-Nancy Cedex\\
FRANCE

\medskip

\noindent Adresse \'electronique : \texttt{deroton@iecn.u-nancy.fr
}


\begin{thebibliography}{99}



\bibitem{D2000} H. Davenport, Multiplicative number theory, 3\up{rd} edition revised by H.L. Montgomery, Springer, 2000.

\bibitem{GG2007} D.A. Goldston et S.M. Gonek, \textit{A note on $S(t)$ and the zeros of the Riemann zeta-function}, Bull. London Math. Soc. {\bf 39} (2007), 482-486.



\bibitem{MM2008} H. Maier et H.L. Montgomery, \textit{The sum of the Möbius function}, \`a para\^{\i}tre dans Journal London Math. Soc. 





\bibitem{S2008} K. Soundararajan, \textit{Partial sums of the Möbius function}, arXiv:0705.0723v2
\bibitem{T2008} G. Tenenbaum, Introduction à la théorie analytique et probabiliste des nombres, 3\up{e} édition, Belin, 2008. 






\end{thebibliography}
\end{document}